# Nonparametric deconvolution problem for dependent sequences*

**Rafał Kulik**

*University of Ottawa*
*e-mail:* rkulik@uottawa.ca

**Abstract:** We consider the nonparametric estimation of the density function of weakly and strongly dependent processes with noisy observations. We show that in the ordinary smooth case the optimal bandwidth choice can be influenced by long range dependence, as opposite to the standard case, when no noise is present. In particular, if the dependence is moderate the bandwidth, the rates of mean-square convergence and, additionally, central limit theorem are the same as in the i.i.d. case. If the dependence is strong enough, then the bandwidth choice is influenced by the strength of dependence, which is different when compared to the non-noisy case. Also, central limit theorem are influenced by the strength of dependence. On the other hand, if the density is supersmooth, then long range dependence has no effect at all on the optimal bandwidth choice.

## 1. Introduction

The nonparametric estimation of the density function for dependent sequences has attracted many researchers in the past. We are not claiming to provide the full overview of this topic, however results can be summarized as follows. In case of weak dependence the results on the (mean square error) optimal bandwidth choice, optimal rates of convergence for the mean square error and central limit theorems for the Parzen-Rosenblatt kernel estimator are exactly the same as in i.i.d. case (see e.g. [3] or [28, Theorem 1]). The situation is a bit more complicated for long-range dependent sequences. Although dependence has no influence on the optimal bandwidth choice, the rates of mean-square convergence may differ according to very strong and moderate dependence. In the latter case they are the same as in the i.i.d. situation. We refer to [5, 8, 13, 14] and [23]. Similarly, if the bandwidth is "small", then central limit theorem for the kernel density estimates is the same as in the i.i.d. case. On the other hand, if the bandwidth is "big" enough, then the long range dependence effect dominates (see [6], [28, Theorem 2]). The similar phenomena occur in random-

*This is an original survey paper

design regression problems. The reader is referred to [22] and [24] for the up-to-date results and references for a kernel and a local linear estimation, respectively.

As for smooth estimators of the distribution function, either short or long range dependence have no influence on the optimal bandwidth choice. The optimal bandwidth is the same as in the i.i.d. case. However, the optimal rates of convergence for the mean square error are always affected by long range dependence (see [8] for more details).

In the present paper we will consider deconvolution problem for dependent sequences. Suppose that we have $n$ observations $Y_1,\ldots,Y_n$ available. We want to estimate the unknown density $f=f_X$ of a random variable $X$, where $Y=X+\epsilon$, with a measurement error $\epsilon$ of a known distribution $F_\epsilon$ and the density $f_\epsilon$. It is assumed that $X$ and $\epsilon$ are independent and that $\{\epsilon,\epsilon_j, j\ge 1\}$ is the i.i.d. sequence.

We will estimate $f(x_0)$ using the classical estimator (cf. [4, 9])

$$\hat{f}_n(x_0)=\frac{1}{nh_n}\sum_{j=1}^{n} g_n\left(\frac{x_0-Y_j}{h_n}\right),$$

where

$$g_n(x_0)=\frac{1}{2\pi}\int_{\mathbf{R}}\exp(-itx_0)\frac{\phi_K(t)}{\phi_\epsilon(t/h_n)}dt.$$

Above, $\phi_\epsilon$ is a characteristic function which corresponds to the density $f_\epsilon$ and $\phi_K(t)=\int_{\mathbf{R}}\exp(itx)K(x)dx$. The mean square error is defined as

$$\mathrm{MSE}(f,h_n):=\mathrm{E}(\hat{f}_n(x_0)-f(x_0))^2.$$

We also study the behavior of the estimator

$$\hat{F}_n(x_0)=\int_{-\infty}^{x_0}\hat{f}_n(u)du.$$

of the distribution function $F$.

In the i.i.d. case the deconvolution problems were studied in [4, 9, 10] and [25] among others. In the latter paper, Fan provided the optimal rates of convergence for $\mathrm{MSE}(f,h_n)$ in both ordinary smooth and supersmooth case. As for weakly dependent case, the previous results have been obtained under various mixing conditions (see [18, 19, 20]) and under association (see [21]). The principal message from the latter papers is that the results (optimal bandwidth, optimal rates, central limit theorem) for weakly dependent sequences are the same as in the i.i.d. case. As for the distribution function, the problem was studied in [9] in the i.i.d case and in [17] in the dependent case.

However, mixing is rather hard to verify and requires additional assumptions. In particular, let $\{Z,Z_i,i\in\mathbb{Z}\}$ be a centered sequence of i.i.d. random variables. Consider the class of stationary linear processes

$$X_j=\sum_{k=0}^{\infty}c_kZ_{j-k},\quad j\ge 1. \tag{1}$$

To obtain a strong mixing for linear processes both regularity of the density of $Z_1$ and some constraints on $c_k$'s are required (see e.g. [7]). On the other hand, association requires that all $c_k$ are positive. To overcome such problems, the martingale based method has been proposed and it works surprisingly well in a variety of problems, not necessary connected with nonparametric estimation (see [15, 16, 26, 27, 28]). Thus, from technical point of view assuming that $c_0 = 1$ and the sequence $c_k$, $k \geq 0$, is summable (referred later as short range dependence, SRD), we will extend Masry's results to moving averages, without referring to mixing or association at all.

However, the more interesting problem is the influence of long range dependence on the deconvolution estimator. To deal with it, we assume that $c_0 = 1$ and $c_k$ is regularly varying with index $-\gamma$, $\gamma \in (1/2, 1)$. This means that $c_k \sim k^{-\gamma} L_0(k)$ as $k \to \infty$, where $L_0$ is slowly varying at infinity. We shall refer to all such models as long range dependent (LRD) linear processes. In particular, if the variance exists, then the covariances $\rho_k := \mathrm{E}X_0X_k$ decay at the hyperbolic rate, $\rho_k = L(k)k^{-(2\gamma-1)}$, where $\lim_{k\to\infty} L(k)/L_0^2(k) = B(2\gamma - 1, 1 - \gamma)$ and $B(\cdot, \cdot)$ is the beta-function. Consequently, the covariances are not summable.

We will show below that in the ordinary smooth case the optimal bandwidth choice for the density problem is influenced by the dependence parameter $\gamma$, as opposite to the optimal bandwidth in the standard (non-noisy) kernel density estimation. In particular, if the dependence is moderate, then the optimal bandwidth and the optimal rates for the density estimation are the same as in the i.i.d. case. If the dependence is very strong, the optimal bandwidth depends on $\gamma$ itself. See Proposition 2.2 and Corollary 2.3. In case of the distribution function, the dependence parameter is always present in the optimal bandwidth and the optimal rates of convergence, as opposite to the non-noisy case (Proposition 2.5).

As for central limit theorem for the density estimator, we have results mimicking CLT for standard kernel density estimators (cf. [28, Theorem 2]): if $h_n$ is small, then CLT is the same as in the i.i.d. case; if $h_n$ is "big", LRD effect starts to dominate. Note that the change from "i.i.d." behavior to LRD behavior occurs in the same way as in the standard kernel estimation, according to $h_n = o\left(\sigma_{n,1}^2/n\right)$ or $\sigma_{n,1}^2/n = o(h_n)$, where

$$\sigma_{n,1}^2 := \mathrm{Var}(\sum_{j=1}^{n} X_j).$$

In the distribution case, we do not have such dichotomous behavior and long-range dependence always influences the limiting behavior.

We note in passing that "small" and "big" bandwidth may have different meanings for different estimation problems. For example, "small" bandwidths are different when estimating a function and its derivative (see [22] for a complete analysis in the regression setting). In the present context and density estimation for error-in-variables models, "small" and "large" bandwidths are the same as for non-noisy case. OF course, this dichotomous behavior is well-known, however, the crucial difference between noisy and non-noisy problem is the optimal bandwidth choice. Note that in the non-noisy setting the opti-mal bandwidth, under appropriate conditions, is *not* influenced by $\gamma$, regardless whether estimation of the function (as mentioned above) or its derivatives is considered. Thus, in errors-in-variables models we have different phenomena than those described in [22].

Another phenomena is that for supersmooth densities the optimal bandwidth choice and the rates for $\mathrm{MSE}(f, h_n)$ are always the same as in the i.i.d. case, irrespectively of the dependence being moderate or very strong. At the first sight this message seems to be optimistic, however, it means that the rate of convergence is so slow that even the very strong dependence cannot worsen it.

## 2. Results

Recall that by SRD assumption we mean that $\sum_{k=0}^{\infty} |c_k| < \infty$. Additionally, we assume that $\sum_{k=0}^{\infty} c_k \neq 0$. By LRD assumption we mean that $c_k \sim k^{-\gamma} L_0(k)$, $\gamma \in (0,1)$.

We assume that $f = f_X$ is twice differentiable with continuous and bounded second order derivatives and $K$ is of the second order. i.e. $\int uK(u)du = 0$ and $0 \neq \int u^2 K(u)du < \infty$. Furthermore, we assume that

$$|\phi_\epsilon(t)| > 0 \tag{2}$$

and that $\phi_\epsilon$, $\phi_K$ are twice differentiable with continuous and bounded derivatives. These assumptions are standard in the i.i.d. situation for both ordinary smooth and supersmooth case.

The proofs are based on the following decomposition: Let $\mathcal{F}_i = \sigma(X_j, Z_j, j \leq i)$. Write

$$\begin{aligned} &\sum_{j=1}^{n} \left( g_n\left(\frac{x_0 - Y_j}{h_n}\right) - \mathrm{E} g_n\left(\frac{x_0 - Y_j}{h_n}\right)\right) \qquad (3) \\ &= \sum_{j=1}^{n} \left( g_n\left(\frac{x_0 - Y_j}{h_n}\right) - \mathrm{E}\left[ g_n\left(\frac{x_0 - Y_j}{h_n}\right) |\mathcal{F}_{j-1}\right]\right) \\ &\quad + \sum_{j=1}^{n} \left( \mathrm{E}\left[ g_n\left(\frac{x_0 - Y_j}{h_n}\right) |\mathcal{F}_{j-1}\right] - \mathrm{E} g_n\left(\frac{x_0 - Y_j}{h_n}\right)\right) =: m_n(x_0) + l_n(x_0). \end{aligned}$$

Note that $\{m_n(x_0), \mathcal{F}_n, n \geq 1\}$ is a martingale. We call $l_n(x_0)$ the *differentiable part*. The similar decomposition is also valid in the distribution case.

### *2.1. Ordinary smooth densities*

Throughout this section, we consider the ordinary smooth case, i.e.

$$|t^{-\beta}||\phi_\epsilon(t)| \to |B_1| > 0, \qquad \beta \geq 1. \tag{4}$$

Furthermore, assume that

$$\delta_{\beta,1} \int |u|^{\beta-2} |\phi_K(u)| du + \int |u|^{\beta-1} |\phi_K'(u)| du + \int |u|^{\beta} |\phi_K^{(2)}(u)| du < \infty, \tag{5}$$

and

$$D_1 := \frac{1}{2\pi|B_1|^2}\int |t|^{2\beta}\phi_K^2(t)dt < \infty. \tag{6}$$

To deal with the LRD case we will impose a stronger condition than (5): with $\beta > 1$,

$$\sum_{j=0}^{3}\int |u|^{\beta-j}|\phi_K^{(j)}(u)| < \infty. \tag{7}$$

We will consider some technical assumptions on the densities:

(A) $f_Z, f_\epsilon$, the densities of $Z$ and $\epsilon$, respectively, are uniformly bounded and Lipschitz continuous.

Note that in this case, $f_{\epsilon+Z}$ - the density $\epsilon+Z$, $f = f_X$ and $f_Y$ - the density of $Y$ are also uniformly bounded and Lipschitz continuous. These conditions are required to handle SRD case.

(B1) $\sum_{r=1}^{2}\int\left|f_{\epsilon+Z}^{(r)}(v)\right|dv < \infty$,

(B2) $\sum_{r=0}^{2}\left(\int\left|f_Z^{(r)}(v)\right|^2 dv + \int\left|f_{\epsilon+Z}^{(r)}(v)\right|^2 dv\right) < \infty$ and $f_Z$, $f_{\epsilon+Z}$ are twice differentiable with continuous and bounded derivatives.

These conditions are required for LRD case (see Appendix for more discussion).

First, we provide the asymptotic expansion of the mean square error.

**Proposition 2.1.** *Assume (4), (5), (6), (A) . Under the SRD assumption,*

$$\begin{aligned}\mathrm{MSE}(f,h_n) \quad &= \quad \left(\frac{1}{2}f^{(2)}(x_0)\int u^2K(u)du\right)^2 h_n^4 + D_1 f_Y(x_0)n^{-1}h_n^{-(2\beta+1)} \\ &\quad + o\left(h_n^4 + n^{-1}h_n^{-(2\beta+1)}\right).\end{aligned}$$

**Proposition 2.2.** *Assume (4), (6), (7), (B1), (B2) and* $\mathrm{E}Z_1^4 < \infty$*. Under the LRD assumption,*

$$\begin{aligned}\mathrm{MSE}(f,h_n) = &\left(\frac{1}{2}f^{(2)}(x_0)\int u^2K(u)du\right)^2 h_n^4 + D_1 f_Y(x_0)n^{-1}h_n^{-(2\beta+1)} + \\ &+ (f_Y'(x_0))^2 n^{-2}\sigma_{n,1}^2 h_n^{-2\beta} + o\left(h_n^4 + n^{-1}h_n^{-(2\beta+1)} + n^{-2}\sigma_{n,1}^2 h_n^{-2\beta}\right).\end{aligned} \tag{8}$$

If $h_n = o(n/\sigma_{n,1}^2)$ (in particular, $h_n = o(n^{-(2(1-\gamma))})$), then the optimal choice is the same as in the i.i.d. case: $h_n = Cn^{-1/(5+2\beta)}$ (here and in the sequel, $C$ is a generic constant, which does not depend on $n$). Consequently, $\mathrm{MSE}(f,h_n) \sim Cn^{-1/(5+2\beta)}$. (We must check that such the choice is permitted, i.e. to check if $n^{-1/(5+2\beta)} = o(n^{-(2(1-\gamma))})$. This is equivalent to $1 > \gamma > \frac{4(2+\beta)+1}{4(2+\beta)+2} > \frac{1}{2}$).

If $n/\sigma_{n,1}^2 = o(h_n)$ (in particular, $n^{-2(1-\gamma)} = o(h_n)$), then we find the optimal bandwidth as $h_n \sim Cn^{-(2\gamma-1)/2(2+\beta)}$. Then, to assure that $n^{-2(1-\gamma)} = o(n^{-(2\gamma-1)/2(2+\beta)})$ we assume that $\gamma < \frac{4(2+\beta)+1}{4(2+\beta)+2}$. From Proposition 2.2 we obtain the following corollary.

**Corollary 2.3.** *Assume that* $\frac{1}{2} < \frac{4(2+\beta)+1}{4(2+\beta)+2} < \gamma < 1$. *Then, by choosing* $h_n = Cn^{-1/(5+2\beta)}$ *we obtain* $\mathrm{MSE}(f, h_n) \sim Cn^{-4/(5+2\beta)}$.

*Assume that* $\frac{1}{2} < \gamma < \frac{4(2+\beta)+1}{4(2+\beta)+2} < 1$. *With* $h_n \sim Cn^{-(2\gamma-1)/2(2+\beta)}$ *we obtain* $\mathrm{MSE}(f, h_n) \sim Cn^{-2(2\gamma-1)/(2+\beta)}$.

**Remark 2.4.** The result of Proposition 2.1 extends the previous ones for $\rho-$ and $\alpha-$mixing (see [19, Lemma 2.1b]) or associated sequences ([21]). In principle, it says that the optimal bandwidth, the rate of convergence of $\mathrm{Var}\hat{f}_n(x_0)$ (and, consequently, of $\mathrm{MSE}(f, h_n)$) for weakly dependent sequences are the same as in the i.i.d. case. Thus is also true for LRD sequences with moderate dependence ($\gamma$ close to 1). On the other hand, if the dependence is very strong, then the bandwidth and the rate of convergence may depend on $\gamma$.

As for the distribution estimator we have the following result.

**Proposition 2.5.** *Assume* (4), (6), (7), (B1), (B2) *and* $\mathrm{E}Z_1^4 < \infty$. *Under the LRD assumption,*

$$\mathrm{MSE}(F, h_n) = \left(\frac{1}{2} f'(x_0) \int u^2 K(u) du\right)^2 h_n^4 + (f_Y(x_0))^2 \frac{\sigma_{n,1}^2}{n^2 h_n^{2\beta}} + o\left(\frac{\sigma_{n,1}^2}{n^2 h_n^{2\beta}} + h_n^4\right).$$

We can see that the optimal bandwidth is $h_n \sim C\left(\sigma_{n,1}^2/n^2\right)^{1/2(\beta+2)}$ and optimal mean square error is of the order $\left(\sigma_{n,1}^2/n^2\right)^{2/(\beta+2)}$. Under weak dependence the optimal bandwidth and the optimal mean square error are proportional to $n^{-1/2(\beta+2)}$ and $n^{-2/(\beta+2)}$, respectively. Consequently, in case of the distribution function the optimal bandwidth and the rates change as soon as we cross the boundary between short- and long range dependence.

As for CLT we have the following results.

**Theorem 2.6.** *Suppose that* $nh_n \to \infty$ *and let* $\sigma^2(x_0) = D_1 f_Y(x_0)$. *Under conditions of Proposition 2.1 we have*

$$n^{1/2} h_n^{\beta+1/2} \left(\hat{f}_n(x_0) - \mathrm{E}\hat{f}_n(x_0)\right) \xrightarrow{\mathrm{d}} N(0, \sigma^2(x_0)).$$

**Theorem 2.7.** *Suppose that* $nh_n \to \infty$ *and let* $\sigma^2(x_0) = D_1 f_Y(x_0)$. *Under conditions of Proposition 2.2 we have*

$$n^{1/2} h_n^{\beta+1/2} \left(\hat{f}_n(x_0) - \mathrm{E}\hat{f}_n(x_0)\right) \xrightarrow{\mathrm{d}} N(0, \sigma^2(x_0))$$

*if* $h_n = o(n/\sigma_{n,1}^2)$, *and*

$$\frac{nh_n^\beta}{\sigma_{n,1}} \left(\hat{f}_n(x_0) - \mathrm{E}\hat{f}_n(x_0)\right) \xrightarrow{\mathrm{d}} N(0, (f_Y'(x_0))^2)$$

*if* $n/\sigma_{n,1}^2 = o(h_n)$. *Under the conditions of Proposition 2.2 we have in either case*

$$\frac{nh_n^\beta}{\sigma_{n,1}} \left(\hat{F}_n(x_0) - \mathrm{E}\hat{F}_n(x_0)\right) \xrightarrow{\mathrm{d}} N(0, (f_Y(x_0))^2).$$

**Remark 2.8.** Theorem 2.6 extends results of [10] and [20]. The result of Theorem 2.7 should be compared with Theorem 2 in [28]. Note that the change from SRD to LRD behavior occurs in the same way as in the standard kernel density case, i.e. by crossing the boundary $h_n \sim n/\sigma_{n,1}^2$. Theorem 2.7 describes the dichotomous behavior of $\hat{f}_n(x_0)$. If $f_Y'(x_0) = 0$, then we may establish trichotomous behavior along the lines of Theorem 3 in [28].

**Remark 2.9.** We shall comment on the assumption $\mathrm{E}Z_1^4 < \infty$. This is necessary for us to use Wu [26] result for empirical processes (see Lemma C below). Instead, we can use Giraitis and Surgailis [12] assumption $\mathrm{E}|Z_1|^{2+\delta} < \infty$ together with additional condition on $f_Z$ (See also Section 2.2). However, it does not solve completely the problem in case $\mathrm{E}Z_1^2 < \infty$.

**Remark 2.10.** It would be desirable to extend the results of, especially, Proposition 2.2 and Theorem 2.7 to the multivariate setting. However, it does not seem to be feasible when using the martingale approximation approach as in the current paper.

**Remark 2.11.** We do not provide CLT for $\hat{F}_n(x_0)$ in the weakly dependent case. The martingale method we use here is based on fact that in the density case the differentiable part is negligible compared to the martingale part, provided that SRD conditions hold (compare (18) with (20)). However, in the distribution case if SRD assumptions are fulfilled, then the martingale part and the differentiable part are of the same order and the method does not apply. We also note that the problem is symmetric in $X$ and $\epsilon$, i.e. instead of assuming that $X_j$ are dependent and $\epsilon_j$ are i.i.d., we may assume that $X_j$ are i.i.d. and $\epsilon_j$ are dependent. What is important in our results is the dependence structure of $Y_j$'s. In [17] it is assumed that $X_j$ is mixing and it is claimed that $\mathrm{Var}\hat{F}_n(x_0)$ has different behavior according to $\epsilon_j$ being dependent or i.i.d. Note, however, that their proof of Lemma 3.2(i) is invalid.

To obtain confidence interval for $\hat{f}_n(x_0)$ we choose appropriate bandwidth to make sure that the variance of the estimator dominates the bias term. In particular, in the LRD case it reads as follows.

**Corollary 2.12.** *Assume that* $h_n = o(n/\sigma_{n,1}^2)$ *and* $h_n = o(n^{-1/(5+2\beta)})$ *(which ensures that* $\frac{1}{2} < \frac{4(2+\beta)+1}{4(2+\beta)+2} < \gamma < 1$*). Then*

$$n^{1/2} h_n^{\beta+1/2} \left( \hat{f}_n(x_0) - f(x_0) \right) \stackrel{\mathrm{d}}{\to} N(0, \sigma^2(x_0)). \tag{9}$$

**Corollary 2.13.** *Assume that* $n/\sigma_{n,1}^2 = o(h_n)$ *and* $h_n = o(n^{-(2\gamma-1)/2(2+\beta)})$ *(which ensures* $\frac{1}{2} < \gamma < \frac{4(2+\beta)+1}{4(2+\beta)+2} < 1$*). Then*

$$\frac{n h_n^{\beta}}{\sigma_{n,1}} \left( \hat{f}_n(x_0) - f(x_0) \right) \stackrel{\mathrm{d}}{\to} N(0, (f_Y'(x_0))^2).$$

To apply Corollary 2.12 in a practical situation one has to estimate $f_Y(x_0)$. If $h_n = o(n/\sigma^2_{n,1})$ and $h_n = o(n^{-1/(5+2\beta)})$, then we can estimate $f_Y(x_0)$ by using

$$\hat{f}_{Y,n}(x_0) = \frac{1}{nb_n}\sum_{j=1}^{n} K\left(\frac{x_0 - Y_j}{b_n}\right).$$

Note that the kernel density bandwidth $b_n$ need not to be the same as $h_n$, as assumed in [20]. In fact, the optimal mean-square error choice is $b_n \sim Cn^{-1/5}$. Further, if $\gamma > 9/10$ (which is ensured by taking $\frac{1}{2} < \frac{4(2+\beta)+1}{4(2+\beta)+2} < \gamma < 1$), then both the mean square error and the variance of the kernel density estimator behave like $Cn^{-4/5}$. On the other hand, from Corollary 2.3, the variance of the deconvolution estimator behaves like $n^{-1}h_n^{-(2\beta-1)}$. Thus, under the constrain $h_n = o(n^{-1/(5+2\beta)})$ the variance of the deconvolution estimator of $f_X(x_0)$ dominates the variance of the kernel density estimator of $f_Y(x_0)$. Consequently, we may build confidence interval for $\hat{f}_n(x_0)$ by replacing $f_Y(x_0)$ with its kernel estimator in (9).

**Remark 2.14.** As suggested by the Referee, assume that

$$X_j = (1-B)^{-\delta_0}\phi^{-1}(B)\psi(B)Z_j,$$

where as before $Z_j$ is i.i.d., $B$ is the backshift operator and $\phi$, $\psi$ are polynomials. If $\delta_0 \in (0, 1/2)$ it is the particular case of LRD model (1) with the specification $1 - \gamma = \delta_0$. If $\delta_0 \in (-1/2, 0)$ (so that $\gamma \in (1, 3/2)$), then this is the case of *antipersistent* sequences. Under appropriate regularity conditions, it was proven in [1, Theorem 3] that for all $\gamma \in (1/2, 3/2)$

$$\frac{1}{n^{3/2-\gamma}}\sum_{i=1}^{n} X_i \xrightarrow{\mathrm{d}} N(0, v),$$

where $v$ is a finite and positive constant. In view of the above result, it is intuitively clear that the expansion (8) is valid for $\gamma \in (1, 3/2)$ as well. However, following the comment below Proposition 2.2, in such the case we always have $h_n = o(n/\sigma^2_{n,1})$. Consequently, the rates of convergence are the same as in i.i.d. case. This is in contrast to fixed-design regression, where antipersistency may improve the rates of convergence beyond those for i.i.d. case. See [1] and [2] for more details.

### *2.2. Supersmooth densities*

In a supersmooth case we consider the usual assumptions (cf. [20]):

(i) $B_1|t|^{\beta_0}\exp(-a|t|^{\beta}) \le |\phi_\epsilon(t)| \le B_2|t|^{\beta_0}\exp(-a|t|^{\beta})$ as $t \to \infty$ for some $a > 0$, $B_1, B_2, \beta > 0$, $\beta_0 \in \mathbb{R}$.

(ii) $\phi_K$ has a finite support $(-d, d)$.

(iii) $|\phi_K(t)| \le B_3(d-t)^m$ and $\phi_K(t) \ge B_4(d-t)^m$ for $t \in (d-\delta, d)$ and positive constants $\delta, m, B_3, B_4$.

(iv) The real (imaginary) part of $\phi_\epsilon$ is negligible as $t \to \infty$ with respect to the imaginary (real) part.

To deal with LRD linear processes (1), recall that $\rho_k = L(k)k^{-(2\gamma-1)}$ as $k \to \infty$. We assume additionally that for all $x, y$,

$$f_k(x,y) = f_Y(x)f_Y(y) + \rho_k f'_Y(x)f'_Y(y) + h_k(x,y), \tag{10}$$

where $f_k$ is the joint density of $(Y_0, Y_k)$ and

$$|h_k(x,y)| \le |\rho_k|^{1+\delta} h(x)h(y), \tag{11}$$

with $\delta > 0$ and $h$ being an integrable and continuous function.

**Proposition 2.15.** *Assume* (i)–(iv) *and that* $f'_Y$ *is continuous and integrable. Under the LRD assumption and* (10) *we have*

$$\mathrm{MSE}(f, h_n) = O((\ln n)^{-2/\beta})$$

*by choosing* $h_n = d\left(\frac{2a}{(1-\theta)\ln n}\right)^{1/\beta}$, $\theta \in (2-2\gamma, 1)$.

The result of Proposition 2.15 means that in the supersmooth case long range dependence has no influence on the optimal bandwidth choice and the optimal rates for $\mathrm{MSE}(f, h_n)$. They are the same as in the i.i.d. and weakly dependent situation situation (cf. [9, 19]).

**Remark 2.16.** Let $f_{k,X}$ ne the joint density of $(X_0, X_k)$, where $\{X_j, i \ge 1\}$ is the LRD linear process $X_j = \sum_{k=0}^{\infty} c_k Z_{j-k}$. Then

$$f_{k,X}(x,y) = f_X(x)f_X(y) + \rho_k f'_X(x)f'_X(y) + \tilde{h}_k(x,y) \tag{12}$$

with $\tilde{h}$ satisfying (11), provided appropriate smoothness condition on $\phi_Z$. We refer to [8, 12] or [23] for more details. Consequently, having established (12), it is easy to verify (10).

**Remark 2.17.** Note that the martingale approximation method used in the ordinary smooth case requires the precise information about $||g_n||_1$, in particular, its finiteness. It is not feasible in the supersmooth case. Instead, we additionally assume (10). We could have worked with this assumption in the ordinary smooth case and obtain the results for $\mathrm{MSE}(f, h_n)$. However, using linear structure and the martingale approximation method we can obtain at the same time $\mathrm{MSE}(f, h_n)$ and the central limit theorem.

## 3. Proofs

Since $f_X$ is twice differentiable with continuous and bounded second order derivatives and $K$ is of the second order, we obtain (see [18])

$$\mathrm{bias}(\hat{f}_n(x_0)) \sim h_n^2 \frac{1}{2} f_X^{(2)}(x_0) \int u^2 K(u)du. \tag{13}$$

### 3.1. Ordinary smooth case

We have

$$\begin{aligned} MSE(h_n) &= \mathrm{Var}\hat{f}_n(x_0) + (\mathrm{bias}(\hat{f}_n(x_0)))^2 \\ &= \frac{1}{(nh_n)^2}\left(\mathrm{E}m_n^2(x_0) + \mathrm{E}l_n^2(x_0) + 2\mathrm{E}m_n(x_0)l_n(x_0)\right) + (\mathrm{bias}(\hat{f}_n(x_0)))^2. \end{aligned}$$

Since $m_n$ is a martingale,

$$\mathrm{E}m_n^2(x_0) = n\mathrm{E}\left(g_n\left(\frac{x_0 - Y_1}{h_n}\right) - \mathrm{E}\left[g_n\left(\frac{x_0 - Y_1}{h_n}\right)|\mathcal{F}_0\right]\right)^2.$$

Let $\zeta_i = g_n\left(\frac{x_0-Y_i}{h_n}\right)$, $\xi_i = \zeta_i - \mathrm{E}(\zeta_i|\mathcal{F}_{i-1})$. Under (2), (4) and (5) we have (cf. [18, Lemma 3])

$$||g_n||_1 = O(h_n^{-\beta}). \tag{14}$$

Also, if additionally (6) holds, then

$$\mathrm{E}[\zeta_i^2] = \int g_n^2\left(\frac{x_0 - u}{h_n}\right) f_Y(u)du \sim D_1 f_Y(x_0)h_n^{1-2\beta} \tag{15}$$

as $n \to \infty$, see [18, Lemma 4].

Let $X_{j,j-1} = \sum_{k=1}^{\infty} c_k Z_{j-k} = X_j - Z_j$. Then, by (14),

$$\begin{aligned} \mathrm{E}(\mathrm{E}[\zeta_i|\mathcal{F}_{i-1}])^2 &= \mathrm{E}\left(\mathrm{E}\left[g_n\left(\frac{x_0 - Y_1}{h_n}\right)|\mathcal{F}_0\right]\right)^2 \\ &= \mathrm{E}\left(h_n\int g_n(x_0 - (u + X_{1,0}))f_{\epsilon+Z}(u)du\right)^2 \\ &\le O(h_n^2)\mathrm{E}\left(\int |g_n(x_0 - (u + X_{1,0}))|du\right)^2 \\ &= O(h_n^2)\mathrm{E}\left(\int |g_n(v)|dv\right)^2 = O(h_n^{2-2\beta}). \end{aligned} \tag{16}$$

Thus, by (15), (16) and Cauchy-Schwartz inequality,

$$\begin{aligned} \mathrm{E}\xi_i^2 &= \mathrm{E}\zeta_i^2 + \mathrm{E}(\mathrm{E}[\zeta_i|\mathcal{F}_{i-1}])^2 - 2\mathrm{E}\left[\zeta_i\mathrm{E}[\zeta_i|\mathcal{F}_{i-1}]\right] \\ &= \mathrm{E}\zeta_i^2 + O(h_n^{2(1-\beta)}) + O(h_n^{\frac{1}{2}-\beta}h_n^{1-\beta}) = \mathrm{E}\zeta_i^2 + o(h_n^{1-2\beta}). \end{aligned} \tag{17}$$

Consequently, via (15), $\mathrm{E}\xi_i^2 \sim D_1 f_Y(x_0)h_n^{1-2\beta}$ as $n \to \infty$ (note that $\xi_i$ depends on $n$) and

$$\mathrm{E}m_n^2(x_0) \sim D_1 f_Y(x_0)nh_n^{1-2\beta}. \tag{18}$$

For $r = 0, 1, 2$, let

$$R_n^{(r)}(z) = \sum_{j=1}^{n}(F_{\epsilon+Z}^{(r)}(x_0 - X_{j,j-1} + z) - F_Y^{(r)}(x_0 + z)).$$

Note that $F_Y^{(r)}(y) = \mathrm{E}F_{\epsilon+Z}^{(r)}(y - X_{1,0})$ (see Lemma D). Also, $\mathrm{E}\left[g_n\left(\frac{x_0 - Y_j}{h_n}\right)|\mathcal{F}_{j-1}\right] = h_n \int g_n(v) f_{\epsilon+Z}(x_0 - X_{j,j-1} - h_n v)dv$. Consequently,

$$l_n(x_0) = h_n \int g_n(v) R_n^{(1)}(-h_n v)dv. \tag{19}$$

*Proof of Proposition 2.1.* We claim that under SRD assumption, $\mathrm{Var}l_n(x_0) = O(nh_n^{2-2\beta})$. Indeed, by [28, Lemma 3], $\sup_z ||R_n^{(1)}(z)||_2^2 = O(n)$. From (19) and by Cauchy-Schwartz inequality,

$$\begin{aligned}
\mathrm{E}l_n^2(x_0) &= h_n^2 \mathrm{E}\left(\int g_n(v) R_n^{(1)}(-h_n v)dv\right)^2 \\
&\leq h_n^2 \mathrm{E} \int\int g_n(u) g_n(v) R_n^{(1)}(-h_n u) R_n^{(1)}(-h_n v)dudv \\
&\leq \int\int |g_n(u)||g_n(v)| \mathrm{E}\left[|R_n^{(1)}(-h_n u)||R_n^{(1)}(-h_n v)|\right] dudv \\
&\leq O(n)\left(\int |g_n(u)|du\right)^2 = O(nh_n^{2-2\beta}).
\end{aligned} \tag{20}$$

Consequently, comparing (20) with (18) we see that $l_n(x_0)$ is negligible compared to the martingale part $m_n(x_0)$. Also, via Cauchy-Schwartz inequality, the mixed term $l_n(x_0)m_n(x_0)$ is negligible. The result of Proposition 2.1 follows by considering (13) and (18). □

*Proof of Proposition 2.2.* Recall (19). Take Taylor expansion,

$$R_n^{(1)}(-h_n v) = R_n^{(1)}(0) - h_n v R_n^{(2)}(\xi),$$

where $\xi = \xi(v) = \xi(v, h_n)$. Thus,

$$\begin{aligned}
\mathrm{E}l_n^2(x_0) &= h_n^2 \mathrm{E}R_n^{(1)}(0)\left(\int g_n(v)dv\right)^2 + O(h_n^4)\mathrm{E}\left(\int v g_n(v) R_n^{(2)}(\xi(u))du\right)^2 \\
&\quad + O(h_n^3)\int g_n(u)du \int u g_n(u) \mathrm{E}\left[R_n^{(1)}(0) R_n^{(2)}(\xi(u))\right] du.
\end{aligned}$$

From Lemmas B, C, the first term is of order $(f_Y'(x_0))^2 h_n^2 \sigma_{n,1}^2 h_n^{-2\beta}$. On account of Lemmas A, (C) and by Cauchy-Schwartz inequality, the third term is bounded by

$$\begin{aligned}
&O(h_n^4)\int\int u g_n(u) v g_n(v) \mathrm{E}\left[R_n^{(2)}(\xi(u)) R_n^{(2)}(\xi(v))\right] dudv \\
&= O(h_n^4) \sup_z \mathrm{E}\left[\left(R_n^{(2)}(z)\right)^2\right]\left(\int |u||g_n(u)|du\right)^2 \\
&= O(h_n^4 \sigma_{n,1}^2)\left(h_n^{\beta+1}\int |u||g_n(u)|du\right)^2 h_n^{-2(\beta+1)} = o(h_n^{2-2\beta}\sigma_{n,1}^2).
\end{aligned}$$

Similarly, for the second term the bound is (cf. Lemmas A and C)

$$
\begin{aligned}
& O(h_n^3 h_n^{-\beta})\left(\mathrm{E}\left(R_n^{(1)}(0)\right)^2\right)^{1/2}\left(\sup_z \mathrm{E}\left[\left(R_n^{(2)}(z)\right)^2\right]\right)^{1/2}\int |ug_n(u)|du \\
= \quad & O(h_n^3 h_n^{-\beta}\sigma_{n,1}^2)h_n^{-(\beta+1)}o(1) = o(h_n^{2-2\beta}\sigma_{n,1}^2).
\end{aligned}
$$

We conclude

$$
\mathrm{E}l_n^2(x_0) = (f_Y'(x_0))^2\sigma_{n,1}^2 h_n^{2-2\beta} + o(\sigma_{n,1}^2 h_n^{2-2\beta}). \tag{21}
$$

Combining (13), (18) and (21) we obtain (8). □

*Proof of Proposition 2.5.* We sketch it briefly, since it is similar to the previous one.

Let $G_n(x_0) = \int_{-\infty}^{x_0} g_n(u)du$. Then $\hat{F}_n(x_0) = \frac{1}{n}\sum_{j=1}^n G_n\left(\frac{x_0-Y_j}{h_n}\right)$ and

$$
\hat{F}_n(x_0) - \mathrm{E}\hat{F}_n(x_0) = \frac{1}{n}\sum_{j=1}^n\left(G_n\left(\frac{x_0-Y_j}{h_n}\right) - \mathrm{E}\left[G_n\left(\frac{x_0-Y_j}{h_n}\right)\right]\right)
$$

Similarly to (3),

$$
\begin{aligned}
&\sum_{j=1}^n\left(G_n\left(\frac{x_0-Y_j}{h_n}\right) - \mathrm{E}G_n\left(\frac{x_0-Y_j}{h_n}\right)\right) \\
&\quad = \sum_{j=1}^n\left(G_n\left(\frac{x_0-Y_j}{h_n}\right) - \mathrm{E}\left[G_n\left(\frac{x_0-Y_j}{h_n}\right)|\mathcal{F}_{i-1}\right]\right) \\
&\qquad + \sum_{j=1}^n\left(\mathrm{E}\left[G_n\left(\frac{x_0-Y_j}{h_n}\right)|\mathcal{F}_{i-1}\right] - \mathrm{E}G_n\left(\frac{x_0-Y_j}{h_n}\right)\right) =: M_n(x_0) + L_n(x_0).
\end{aligned}
$$

Then,

$$
\mathrm{E}M_n^2(x_0) = O(nh_n^{-2\beta}). \tag{22}
$$

Further, we have as in (19),

$$
L_n(x_0) = h_n\int G_n(v)R_n^{(0)}(-h_n v)dv = \int g_n(v)R_n^{(0)}(-h_n v)dv.
$$

Taking Taylor expansion $R_n^{(0)}(-h_n v) = R_n^{(0)}(0) - h_n v R_n^{(1)}(\xi)$ we obtain as in the proof of Proposition 2.2,

$$
\mathrm{E}L_n^2(x_0) = (f_Y(x_0))^2\sigma_{n,1}^2 h_n^{-2\beta} + o\left(\sigma_{n,1}^2 h_n^{-2\beta}\right). \tag{23}
$$

Comparing (22) with (23) we can see that the martingale part is of smaller order. Consequently,

$$
\mathrm{Var}\hat{F}_n(x_0) = (f_Y(x_0))^2\frac{\sigma_{n,1}^2}{n^2}h_n^{-2\beta} + o\left(\frac{\sigma_{n,1}^2}{n^2}h_n^{-2\beta}\right).
$$

Since $\mathrm{bias}(\hat{F}_n(x_0)) = \frac{1}{2}f'(x_0)\int u^2K(u)du h_n^2 + o(h_n^2)$ we conclude the result. □

In order to prove CLT, we will use the martingale central limit theorem.

**Lemma 3.1.** *Assume that* $nh_n \to \infty$. *Then*

$$(nh_n^{1-2\beta})^{-1/2} m_n(x_0) \xrightarrow{\mathrm{d}} N(0, \sigma^2(x_0))$$

*Proof.* The proof is similar to that of Lemma 2 in [28].

Since $m_n$ is a martingale it suffices to verify the Lindeberg condition and convergence of conditional variances.

Let $\bar{\zeta}_j = (nh_n^{1-2\beta})^{-1/2}\zeta_i$, $\bar{\xi}_j = (nh_n^{1-2\beta})^{-1/2}\xi_i$. Note that for sufficiently large $n$ we have $h_n^\beta |g_n(v)| \le C$ and the bound does not depend on $v$ nor $n$. As for the Lindeberg condition we have by

$$\begin{aligned} n\mathrm{E}\left[\bar{\xi}_j^2 1_{\{|\bar{\xi}_i|>\varepsilon\}}\right] &\le 4n\mathrm{E}\left[\bar{\zeta}_j^2 1_{\{|\bar{\zeta}_j|>\varepsilon/2\}}\right] \\ &= 4h_n^{2\beta}\int g_n^2(v) f_Y(x_0 - vh_n) 1_{\{|g_n(v)|>(nh_n^{1-2\beta})^{1/2}\varepsilon/2\}} dv \\ &= O(h_n^{2\beta}) \int g_n^2(v) 1_{\{h_n^\beta |g_n(v)|>(nh_n)^{1/2}\varepsilon/2\}} dv. \end{aligned}$$

The set $\{h_n^\beta |g_n(v)| > (nh_n)^{1/2}\varepsilon/2\}$ becomes empty for sufficiently large $n$. Consequently,

$$n\mathrm{E}\left[\bar{\xi}_j^2 1_{\{|\bar{\xi}_j|>\varepsilon\}}\right] \to 0$$

as $n \to \infty$.

Now, we want to show that

$$\sum_{j=1}^n \mathrm{E}\left[\bar{\xi}_j^2 | \mathcal{F}_{j-1}\right] \xrightarrow{p} \sigma^2(x_0). \tag{24}$$

As in (17), we have

$$\mathrm{E}\left[\bar{\xi}_j^2 | \mathcal{F}_{j-1}\right] = \mathrm{E}\left[\bar{\zeta}_j^2 | \mathcal{F}_{j-1}\right] - \mathrm{E}\left(\mathrm{E}[\bar{\zeta}_j | \mathcal{F}_{j-1}]\right)^2 = \mathrm{E}\left[\bar{\zeta}_j^2 | \mathcal{F}_{j-1}\right] + O_P(n^{-1}h_n).$$

Consequently,

$$\sum_{j=1}^n \mathrm{E}\left[\bar{\xi}_j^2 | \mathcal{F}_{j-1}\right] - \sum_{j=1}^n \mathrm{E}\left[\bar{\zeta}_j^2 | \mathcal{F}_{j-1}\right] = o_P(1)$$

and it suffices to prove

$$\sum_{j=1}^n \mathrm{E}\left[\bar{\zeta}_j^2 | \mathcal{F}_{j-1}\right] \xrightarrow{p} \sigma^2(x_0).$$

We have

$$\begin{aligned} &\sum_{j=1}^n \mathrm{E}\left[\bar{\zeta}_j^2 | \mathcal{F}_{j-1}\right] - \sum_{j=1}^n \mathrm{E}\left[\bar{\zeta}_j^2\right] \\ &= h_n^{2\beta} \int g_n^2(v) \frac{1}{n} R_n^{(1)}(-h_n v) dv \\ &\le h_n^{2\beta} \int g_n^2(v) \frac{1}{n} |R_n^{(1)}(-h_n v) - R_n^{(1)}(0)| dv + \frac{1}{n} R_n^{(1)}(0) h_n^{2\beta} \int g_n^2(v) dv. \end{aligned}$$

By (15) and the ergodic theorem, the second part is $o_P(1)$. By Lipschitz continuity of $f_{\epsilon+Z}$ and $f_Y$, the first part is bounded by

$$O(h_n^{2\beta})\int g_n^2(v)\min\{1,|h_nv|\}dv$$

which converges to 0 by the dominated convergence theorem since $h_n^{2\beta}g_n(v)$ is uniformly bounded in $v$ and $n$ and integrable. □

*Proof of Theorem 2.6.* It follows from Lemma 3.1 and (20). □

*Proof of Theorem 2.7.* If $h_n = o(n/\sigma_{n,1})$, then it follows from Lemma 3.1 and (21), since the martingale part dominates the differentiable part.

If $n/\sigma_{n,1}^2 = o(h_n)$, then as in the proof of Proposition 2.2,

$$\sigma_{n,1}^{-1}h_n^{1-\beta}l_n(x_0) = R_n^{(1)}(0) + o_P(1).$$

Consequently, the result follows by Lemma C since the martingale part is negligible. □

### 3.2. Supersmooth case

*Proof of Proposition 2.15.* Let $Z_{n,j} = \frac{1}{h_n}g_n\left(\frac{x_0-Y_j}{h_n}\right)$. Then

$$\begin{aligned}\mathrm{Var}\hat{f}_n(x_0) &= \frac{1}{n}\mathrm{Var}Z_{n,0} + \frac{2}{n}\sum_{j=1}^{n-1}(1-j/n)\mathrm{Cov}(Z_{n,0},Z_{n,j}) \\ &= \frac{1}{n}\mathrm{Var}Z_{n,0} + O\left(\frac{2}{n}\sum_{j=1}^{n-1}(1-j/n)|\rho_j|^{1+\delta}\right)\frac{1}{h_n^2}\left(\int g_n\left(\frac{x_0-u}{h_n}\right)h(u)du\right)^2 \\ &\quad + \frac{2}{n}\sum_{j=1}^{n-1}(1-i/n)\rho_j\frac{1}{h_n^2}\left(\int g_n\left(\frac{x_0-u}{h_n}\right)f_Y'(u)du\right)^2 =: I_1+I_2+I_3.\end{aligned}$$

From [20] we know that

$$l_n^{(1)} := C\frac{1}{n}h_n^{2([m+1]\beta+\beta_0-\frac{1}{2})}\exp(2a(d/h_n)^\beta) \le I_1$$

and

$$I_1 \le C\frac{1}{n}h_n^{2([m+1]\beta+\beta_0-1)}\exp(2a(d/h_n)^\beta)\left(\ln(1/h_n)\right)^{2m} =: U_n^{(1)}.$$

Now, from continuity and integrability of $f_Y'$ we obtain via Lemma 3.1 in [11],

$$\frac{1}{h_n^2}\left(\int g_n\left(\frac{x_0-u}{h_n}\right)f_Y'(u)\right)^2 \ge Ch_n^{2([m+1]\beta+\beta_0)}\exp(2a(d/h_n)^\beta)$$

and

$$\frac{1}{h_n^2}\left(\int g_n\left(\frac{x_0-u}{h_n}\right)f_Y'(u)\right)^2 \le Ch_n^{2([m+1]\beta+\beta_0-1)}\exp(2a(d/h_n)^\beta)\,(\ln(1/h_n))^{2m}. \tag{25}$$

Consequently,

$$l_n^{(3)} := C\frac{\sigma_{n,1}^2}{n^2}h_n^{2([m+1]\beta+\beta_0-\frac{1}{2})}\exp(2a(d/h_n)^\beta) \le I_3$$

and

$$I_3 \le C\frac{\sigma_{n,1}^2}{n^2}h_n^{2([m+1]\beta+\beta_0-1)}\exp(2a(d/h_n)^\beta)\,(\ln(1/h_n))^{2m} =: U_n^{(3)}.$$

Further, as in (25),

$$I_2 \le o\left(\frac{\sigma_{n,1}^2}{n^2}\right)h_n^{2([m+1]\beta+\beta_0-1)}\exp(2a(d/h_n)^\beta)\,(\ln(1/h_n))^{2m} =: U_n^{(2)}.$$

To assure that $I_1 \to 0$ as $n \to \infty$ we choose $h_n = d\left(\frac{2a}{(1-\theta)\ln n}\right)^{1/\beta}$, $\theta \in (0,1)$. Now, $U_n^{(1)} = o(l_n^{(3)})$ and $U_n^{(2)} = o(l_n^{(3)})$ as long as $h_n \sim C(\ln n)^{-\kappa}$, $\kappa > 0$. Consequently, with our choice of $h_n$ the third part $I_3$ dominates both $I_1$ and $I_2$. The upper bound for the third part is

$$O\left(\frac{\sigma_{n,1}^2}{n^{1+\theta}}\right) = O(n^{-\kappa}),$$

$\kappa > 0$ as long as $0 < 2-2\gamma < \theta < 1$. Consequently, via (13), the bias term dominates and the mean square rate of convergence is of the order $(\ln n)^{-2/\beta}$. □

## Appendix A

**Lemma A.** *Assume* (2), (4), (7). *Then* $ug_n(u) \in L_1(\mathbb{R})$, $\int |ug_n(u)|du = O(h_n^{-2})$ *and consequently* $h_n^{\beta+1}\int |ug_n(u)|du = o(1)$ *for* $\beta > 1$.

*Proof.* Integrate by parts three times to obtain

$$|u^3 g_n(u)| < \frac{1}{2\pi}\int \left|\left(\frac{\phi_K(t)}{\phi_\epsilon(t/h_n)}\right)^{(3)}\right| dt.$$

Consequently, if we show that the right-hand side is bounded by $Ch_n^{-2}$ we will prove that $|ug_n(u)| = O(|u|^{-2})$ (the bound depends on $h_n$) and hence $|ug_n(u)|$ is integrable.

We have

$$\left(\frac{\phi_K(t)}{\phi_\epsilon(t/h_n)}\right)^{(3)} = \frac{\phi_K^{(3)}(t)}{\phi_\epsilon(t/h_n)} - \frac{3}{h_n}\frac{\phi_K^{(2)}(t)\phi'_\epsilon(t/h_n)}{\phi_\epsilon^2(t/h_n)} + \frac{5}{h_n^2}\frac{\phi'_K(t)(\phi'_\epsilon(t/h_n))^2}{\phi_\epsilon^3(t/h_n)}$$
$$- \frac{1}{h_n^2}\frac{\phi_\epsilon^{(2)}(t/h_n)}{\phi_\epsilon^2(t/h_n)}[2\phi_K(t) + \phi'_K(t)] + \frac{4}{h_n^3}\frac{\phi_K(t)\phi_\epsilon^{(2)}(t/h_n)\phi'_\epsilon(t/h_n)}{\phi_\epsilon^3(t/h_n)}$$
$$- \frac{1}{h_n^3}\frac{\phi_K(t)}{\phi_\epsilon^2(t/h_n)}\left[\phi_\epsilon^{(3)}(t/h_n) + \frac{(\phi'_\epsilon(t/h_n))^3}{\phi_\epsilon^2(t/h_n)}\right].$$

Taking integral of each term separately on $\{|t| < Mh_n\}$ and $\{|t| \geq Mh_n\}$, using (2), (4) and boundness of derivatives we obtain for the terms involving $h_n^3$:

$$\frac{4}{h_n^3}\int_{\{|t|<Mh_n\}}\left|\frac{\phi_K(t)\phi_\epsilon^{(2)}(t/h_n)\phi'_\epsilon(t/h_n)}{\phi_\epsilon^3(t/h_n)}\right|dt$$
$$+ \frac{1}{h_n^3}\int_{\{|t|<Mh_n\}}\left|\frac{\phi_K(t)}{\phi_\epsilon^2(t/h_n)}\left[\phi_\epsilon^{(3)}(t/h_n) + \frac{(\phi'_\epsilon(t/h_n))^3}{\phi_\epsilon^2(t/h_n)}\right]\right|dt = O(h_n^{-2})$$

On $\{|t| \geq Mh_n\}$ we utilize condition (7) and the form (2) of $\phi_\epsilon$ and the corresponding behavior of its derivatives. □

To establish exact asymptotics in the LRD case, we need the precise result on behavior of $\in g_n(u)du$.

**Lemma B.** *Assume* (2), (4), (5). *Then*

$$\lim_{n\to\infty} h_n^\beta \int g_n(u)du = 1.$$

*Proof.* In view of [18, Lemma 3], $g_n \in L_1(\mathbb{R})$. Let $g(u) = \int \exp(itu)\frac{\phi_K(t)}{\phi_\epsilon(t)}dt$. By the inversion formula, $\phi_K(t)/\phi_\epsilon(t) = \int \exp(-itu)g(u)du$. Since $g \in L_1(\mathbb{R})$, taking $t = 0$, we obtain $\int g(u)du = 1$.

On the other hand,

$$\lim_{n\to\infty}\int h_n^\beta g_n(u)du = \int \lim_{n\to\infty} h_n^\beta g_n(u)du$$
$$= \frac{1}{2\pi}\int\int \exp(itu)\lim_{n\to\infty} h_n\frac{\phi_K(t)}{\phi_\epsilon(t/h_n)}dtdu.$$

The change of limit and integrals is permitted since ([18, Lemma 3])

$$h_n^\beta|g_n(u)| = h_n^\beta|g_n(u)|1_{\{|u|<1\}} + h_n^\beta|g_n(u)|1_{\{|u|\geq 1\}}$$
$$\leq Ch_n^\beta 1_{\{|u|<1\}}\int\left|\frac{\phi_K(t)}{\phi_\epsilon(t/h_n)}\right|dt + Ch_n^{-\beta}u^{-2}1_{\{|u|\geq 1\}}$$

and

$$h_n^\beta\left|\frac{\phi_K(t)}{\phi_\epsilon(t/h_n)}\right| = h_n^\beta\left|\frac{\phi_K(t)}{\phi_\epsilon(t/h_n)}\right|1_{\{|t|<Mh_n\}} + h_n^\beta\left|\frac{\phi_K(t)}{\phi_\epsilon(t/h_n)}\right|1_{\{|t|\geq Mh_n\}}$$
$$\leq Ch_n^\beta 1_{\{|t|<Mh_n\}} + Ch_n^\beta|\phi_K(t)||t^{-\beta}|.$$

The upper bounds are integrable as the functions of $u$ and $t$, respectively. □

## Appendix B

**Lemma C.** *Assume* (B1), (B2) *and* $\mathrm{E}Z_1^4 < \infty$. *Then under the LRD assumptions and* $r = 0, 1, 2$ *we have a weak convergence*

$$\sigma_{n,1}^{-1} R_n^{(r)}(z) \Rightarrow f_Y^{(r+1)}(x_0 + z) Z,$$

$$\sup_z \mathrm{E}\left[\left|R_n^{(r)}(z)\right|^2\right] = O(\sigma_{n,1}^2),$$

*and* $\mathrm{E}\big[\big(R_n^{(r)}(0)\big)^2\big] \sim \sigma_{n,1}^2 \big(f_Y^{(r+1)}(x_0)\big)^2$.

*Proof.* Let $r = 0, 1, 2$. Let $G_n$ be the empirical distribution function associated with $X_{1,0}, \ldots, X_{n,n-1}$. Let $G$ be the distribution of $X_{1,0}$. Then

$$R_n^{(r)}(z) = n \int (G_n(u) - G(u))\, f_{\epsilon+Z}^{(r+1)}(x_0 - u + z) du$$

Consequently, under the condition $\int f_{\epsilon+Z}^{(r+1)}(v) dv < \infty$ we have

$$|R_n^{(r)}(z)| \le Cn \sup_u |G_n(u) - G(u)|$$

and the bound is independent of $z$. Now, we apply Theorem 2 in [26] with $p = 0$. Then $\mathrm{E}\left[\sup_u |G_n(u) - G(u)|^2\right] = O(\sigma_{n,1}^2 / n^2)$ and hence $\sup_z \mathrm{E}\big[|R_n^{(r)}(z)|^2\big] = O(\sigma_{n,1}^2)$.

Further, we can apply Theorem 1 in [26] to obtain

$$\mathrm{E}\left[\left|R_n^{(r)}(0) + f_Y^{(r+1)}(x_0) \sum_{j=1}^{n} X_{j,j-1}\right|^2\right] = o(\sigma_{n,1}^2).$$

Consequently,

$$\sigma_{n,1}^2 \mathrm{E}\left[\left(R_n^{(r)}(0)\right)^2\right] \to \left(f_Y^{(r+1)}(x_0)\right)^2.$$

□

**Lemma D.** *Assume that* $\sum_{r=0}^{2} \int \int \left|f_{\epsilon+Z}^{(r)}(v)\right|^2 dv < \infty$ *and* $\mathrm{E}|Y_1|^\kappa < \infty$ *for some* $\kappa > 0$. *Then*

$$F_Y^{(r)}(y) = \mathrm{E}F_{\epsilon+Z}^{(r)}(y - X_{1,0}).$$

*Proof.* It follows from [26, Lemma 6]. □

**Lemma E.** *Let* $k \ge 1$. *Assume that either* $\int |f_\epsilon^{(r)}(v)|^k dv < \infty$ *or* $\int |f_Z^{(r)}(v)|^k dv < \infty$. *Then* $\int |f_{\epsilon+Z}^{(r)}(v)|^k dv < \infty$.

*Proof.* By Fubini's theorem

$$\int |f_{\epsilon+Z}^{(r)}(v)|^k dv = \int |\mathrm{E}f_Z^{(r)}(v - \epsilon)|^k dv \le \mathrm{E} \int |f_Z^{(r)}(v - \epsilon)|^k dv < \infty.$$

□

**Acknowledgement.**

I would like to thank both referees for valuable comments, which improved the paper.

## References

[1] Beran, J. and Feng, Y. (2001). Local polynomial estimation with a FARIMA-GARCH error process. *Bernoulli*, **7**, 733–750. MR1867080

[2] Beran, J. and Feng, Y. (2001). Local polynomial fitting with long-memory, short-memory and antipersistent errors. *Ann. Inst. Statist. Math.*, **54**, 291–311. MR1910174

[3] Bosq, D. (1996). *Nonparametric Statistics for Stochastic Processes.* Lecture Notes in Statistics **110**. Springer, New York. MR1441072

[4] Carroll, R.J. and Hall, P. (1988). Optimal Rates of Convergence for Deconvoling a Density. *J. Amer. Statist. Assoc.* **83**, 1184–1186. MR0997599

[5] Claeskens, G. and Hall, P. (2002). Effect of dependence on stochastic measures of accuracy of density estimators. *Ann. Statist.* **30**, 431–454. MR1902894

[6] Csőrgö, S. and Mielniczuk, J. (1995). Density estimation under long-range dependence. *Ann. Statist.* **23**, 990–999. MR1345210

[7] Doukhan, P. (1984). *Mixing: Properties and Examples.* Lecture Notes in Statisitcs. Springer. MR1312160

[8] Estevez, G. and Vieu, P. (2003). Nonparametric estimation under long memory dependence. *Nonparametirc Statistics* **15**, 535–551. MR2017486

[9] Fan, J. (1991). On the optimal rates of convergence for nonparametric deconvolution problems. *Ann. Statist.* **19**, 1257–1272. MR1126324

[10] Fan, J. (1991). Asymptotic normality for deconvolving kernel density estimators. *Sankya Ser. A* **53**, 97–110. MR1177770

[11] Fan, J. and Masry, E. (1992). Multivariate Regression Estimation with Errors-in-Variables: Asymptotic Normality for Mixing Processes. *J. Mult. Anal.* **43**, 237–271. MR1193614

[12] Giraitis, L. and Surgailis, D. (1999). Central limit theorem for the empirical process of a linear sequence with long memory. *J. Statist. Plann. Inference* **80**, 81–93. MR1713796

[13] Hall, P. and Hart, J.D. (1990). Convergence rates in density estimation for data from infinite-order moving average processes. *Probab. Th. Rel. Fields* **87**, 253–274. MR1080492

[14] Hall, P., Lahiri, S. N. and Truong, Y. K. (1995). On bandwidth choice for density estimation with dependent data. *Ann. Statist.* **23**, 2241–2263. MR1389873

[15] Ho, H.-C. and Hsing, T. (1996). On the asymptotic expansion of the empirical process of long-memory moving averages. *Ann. Statist.* **24**, 992–1024. MR1401834

[16] Ho, H.-C. and Hsing, T. (1997). Limit theorems for functionals of moving averages. *Ann. Probab.* **25**, 1636–1669. MR1487431

[17] Ioannides, D.A. and Papanastassiou, D.P. (2001). Estimating the Distribution Function of a Stationary Process Involving Measurement Errors. *Statist. Inf. for Stoch. Proc.* **4**, 181–198. MR1856173

[18] Masry, E. (1991). Multivariate Probability Density Deconvolution for Stationary Random Processes. *IEEE Trans. Inf. Th.* **37**, 1105–1115. MR1111811

[19] Masry, E. (1993). Strong consistency and rates for deconvolution of multivariate densities of stationary processes. *Stochastic Process. Appl.* **47**, 53–74. MR1232852

[20] Masry, E. (1993). Asymptotic normality for deconvolving estimators of multivariate densities for stationary processes. *J. Mult. Anal.* **44**, 47–68. MR1208469

[21] Masry, E. (2003). Deconvolving Multivariate Kernel Density Estimates From Contaminated Associated Observations. *IEEE Trans. Inf. Th.* **49**, 2941–2952. MR2027570

[22] Masry, E. and Mielniczuk, J. (1999). Local linear regression estimation for time series with long-range dependence. *Stoch. Proc. Appl.* **82**, 173–193. MR1700004

[23] Mielniczuk, J. (1997). On the asymptotic mean integrated squared error of a kernel density estimator for dependent data. *Statist. Probab. Letters* **34**, 53–58. MR1457496

[24] Mielniczuk, J. and Wu, W. B. (2004). On random-design model with dependent errors. *Statistica Sinica*, 1105–1126. MR2126343

[25] Stefanski, L.A. and Carroll, R.J. (1990). Deconvoltuing kernel density estimators. *Statistics* **21**, 169–184. MR1054861

[26] Wu, W.B. (2003). Empirical processes of long-memory sequences. *Bernoulli* **9**, 809–831. MR2047687

[27] Wu, W.B. (2005). On the Bahadur representation of sample quantiles for dependent sequences. *Ann. Statist.* **33**, 1934–1963. MR2166566

[28] Wu, W. B. and Mielniczuk, J. (2002) Kernel density estimation for linear processes. *Ann. Statist.* **30**, 1441–1459. MR1936325